\documentclass{article}
\usepackage{amssymb,latexsym}
\usepackage{amsmath}

\newtheorem{theorem}{Theorem}
\newtheorem{remark}{Remark}
\def\RR{\mathbb{R}}
\def\NN{\mathbb{N}}

\def\PP{\mathbb{P}}
\def\ZZ{\mathbb{Z}}
\begin{document}

\title{\bf On two families of near-best spline quasi-interpolants on non-uniform partitions of the real line}
\date{January 2006}
\author{  D. Barrera, M.J. Ib\'a\~nez, P. Sablonni\`ere, D. Sbibih}
\maketitle


\begin{abstract}
The univariate spline  quasi-interpolants (abbr. QIs) studied in this paper
are approximation operators using B-spline expansions with coefficients which are
linear combinations of discrete or weighted mean values of the function to be approximated.
When working with nonuniform partitions, the main challenge is to find QIs which have 
both good approximation orders and uniform norms which are bounded independently of the given
partition. Near-best QIs are obtained by minimizing an upper
bound of the infinity norm of QIs depending on a certain number of free
parameters, thus reducing this norm. This paper is devoted to the study of two families of near-best QIs of approximation order 3.
\end{abstract}


\section{Introduction}

A spline quasi-interpolant (abbr. QI) of $f$ has the general form
$$
Qf=\sum_{\alpha\in A} \lambda_{\alpha}(f) B_{\alpha}
$$
where $\{B_{\alpha},\alpha\in A\}$ is a family of B-splines forming a partition of unity and 
$\{\lambda_{\alpha},\alpha\in A\}$
is a family of linear functionals which are local in the sense that they only use values of $f$
in some neighbourhood of  $\Sigma_{\alpha}=$supp$(B_{\alpha})$. The main interest of QIs is that they
provide good approximants of functions without solving any linear system of equations.
In this paper, we want to study  the following types of QIs:\\

{\sl Discrete Quasi-Interpolants} (abbr. dQIs) : the linear functionals are {\sl linear combinations of
values} of $f$ at some points in a neighbourhood of $\Sigma_{\alpha}$ (see e.g. \cite{bis}-\cite{biss2}, \cite{dbf}-\cite{dvl},\cite{mj},\cite{lls}\cite{ls}\cite{sab3}\cite{sab4}).

{\sl Integral Quasi-Interpolants} (abbr. iQIs) : the linear functionals are
{\sl linear combinations of weighted mean values} of $f$ in some neighbourhood of
$\Sigma_{\alpha}$ (see e.g.\cite{biss1}-\cite{db2},\cite{ls},\cite{sab3}-\cite{ss}).\\

More specifically,  we study  QIs  that we call {\sl Near-Best Quasi-Interpolants} (abbr. NB QIs) which are defined as follows:\\

{\bf 1) Near-Best dQIs:} assume that $\lambda_{\alpha}(f)=\sum_{\beta\in F_{\alpha}} a_\alpha(\beta)
f(x_\beta)$ where the finite set of points $\{x_\beta, \beta\in F_{\alpha}\}$ lies in some
neighbourhood of $\Sigma_{\alpha}$. Then it is clear that, for 
$\Vert f \Vert_{\infty}\le 1$ and $\alpha\in A$, 
$\vert \lambda_{\alpha}(f) \vert \le \Vert a_\alpha \Vert_1$,
where $a_\alpha$ is the vector with components $a_\alpha(\beta)$, from which we deduce
immediately
$$
\Vert Q \Vert_{\infty} \le \sum_{\alpha\in A} \vert \lambda_{\alpha}(f) \vert B_{\alpha}\le
\max_{\alpha\in A}\vert \lambda_{\alpha}(f) \vert\le \max_{\alpha\in A}\Vert a_\alpha
\Vert_1=\nu_1(Q).
$$
Now, assuming that $n=$card$(F_\alpha)$ for all $\alpha$, we can try to find
$a_\alpha^*\in\RR^n$ solution of  the minimization problem 
$$
\Vert a_\alpha^*\Vert_1=\min\{\Vert a_\alpha\Vert_1;\; a_\alpha\in\RR^n, \;
V_\alpha a_\alpha=b_\alpha\}
$$
where the linear constraints express that $Q$ is exact on some subspace of polynomials. Thus, we finally obtain 
$$
\Vert Q \Vert_{\infty} \le \nu_1^*(Q)=\max_{\alpha\in A}\Vert a_\alpha^* \Vert_1.
$$
{\bf 2) Near-Best iQIs:} assume that $\lambda_{\alpha}(f)=\sum_{\beta\in F_{\alpha}} a_\alpha(\beta)
\int_{\Sigma_\beta} M_{\beta}(t)f(t)dt$, where the B-splines $M_{\beta}$ are normalized by $\int
M_{\beta}=1$. Note that the B-spline $M_{\beta}$ can be different from $B_{\alpha}$. Once
again, for
$\Vert f
\Vert_{\infty}\le 1$, we have
$$
\vert \lambda_{\alpha}(f)\vert\le \sum_{\beta\in F_{\alpha}} \vert a_\alpha(\beta)\vert
\vert\int_{\Sigma_\beta} M_{\beta}(t)f(t)dt \vert \le \sum_{\beta\in F_{\alpha}}
\vert a_\alpha(\beta)\vert=\Vert a_\alpha \Vert_1
$$
whence, as we obtained above for dQIs,
$$
\Vert Q \Vert_{\infty} \le \max_{\alpha\in A} \Vert a_\alpha \Vert_1=\nu_1^*(Q).
$$
As emphasized by de Boor (see e.g. \cite{db1}, chapter XII), a QI defined on non uniform partitions has
to be {\sl uniformly bounded independently of the partition} in order to be
interesting for applications. Therefore, the aim of this
paper is to define some families of discrete and integral QIs satisfying this property and having the {\sl smallest possible norm}. As in general it is difficult to minimize the true norm of the operator, {\sl we have
chosen to solve the minimization problems} defined above. \\

The paper extends some results of \cite{bis} \cite{mj}, and is organized as follows. We first recall some
"classical" QIs of various types and we verify that they are uniformly bounded. Then we define and
study several families of discrete and integral QIs, depending on a finite number of parameters, for
which we can find $\nu_1^*(Q)$. We show that this problem has always a solution (in general non
unique). 
Of particular interest are the results  of theorems 1,4,6 and 8 where we show that some families of dQIs and iQIs are uniformly bounded independently of the partition. By imposing more constraints on the non-uniform partitions, we can also prove that some families of QIs are near-best (theorems 5 and 9).
(A parallel study of spline QIs is done in \cite{biss1} for uniform partitions of the real line and in \cite{biss2}  for a uniform triangulation of the plane). In all  cases, the QIs that we study are only exact on $\PP_2$, i.e.  their approximation order is $3$. It seems surprisingly difficult to construct QIs which are both uniformly bounded  independently of the partition and exact on $\PP_d$ for $d\ge 3$. In the last section 11, we consider an example of QI which is exact on $\PP_3$ (i.e. of approximation order $4$) and uniformly bounded: however  the bound depends on the {\sl maximal mesh ratio} of the partition. Such operators seem also useful for applications and it would be interesting to study near-best operators of this type, thus allowing a reduction of the upper bound of the norm.


\section{Notations} 

We shall use classical B-splines of degree $m$ on a bounded interval
$I=[a,b]$ or on $I=\RR$. For the sake of simplicity, in the case $I=\RR$, we
take a strictly increasing sequence of knots $T=\{t_i, i\in\ZZ\}$ satisfying  $\vert t_i\vert \to+\infty$
 as $\vert i \vert \to+\infty $.
 In the  case $I=[a,b]$, we take
the usual sequence $T$ of knots defined by (see e.g. \cite{db1},\cite{dvl}) :
$$
t_{-m}=\cdots=t_0=a<t_1<t_2<\cdots<t_{n-1}<b=t_n=\cdots=t_{n+m}.
$$
For $J=\{0,\ldots,n+m-1\}$,
the family of B-splines $\{B_j,j\in J\}$,
with support $\Sigma_j=[t_{j-m},t_{j+1}]$ is a basis of
the space
$S_m(I,T)$ of splines of degree $m$ on the interval $I$ endowed with the
partition $T$. 
These B-splines form a partition of unity, i.e. $\sum_{j\in J} B_j=1$. We set $h_i=t_i-t_{i-1}$
for all indices $i$.

Let $\NN_m=\{1,\ldots,m\}$ and $T_j=\{t_{j-r+1}, r\in \NN_m\}$: 
we recall that the {\sl elementary symmetric functions}
$\sigma_l(T)$  of the $m$ variables in $T_j$  are defined by $\sigma_0(T_j)=0$ and for
$1\le l\le m$, by
$$
\sigma_l(T_j)=\sum_{1\le r_1< r_2<\ldots<r_l\le m}t_{j+1-r_1}t_{j+1-r_2}\ldots t_{j+1-r_l}.
$$
Denoting $C_m^l=\frac{m!}{l!(m-l)!}$ the binomial coefficients, then, for $0\le l\le m$,  the monomials 
$e_l(x)=x^l$ can be written
$e_l=\sum_{i\in J}\theta_i^{(l)} B_i$, with $\theta_i^{(l)}=\sigma_l(T_i)/C_m^l$.
This is a direct consequence of Marsden's identity (\cite{db1}, chapter IX). 
In particular, the Greville points have abscissas $\theta_i=\theta_i^{(1)}.$\\
Similarly, we define the {\sl extended symmetric functions}
$\bar \sigma_l(T_j)$ by $\bar\sigma_0(T_j)=1$ and, for $1\le l\le m$, 
$$
\bar \sigma_l(T_j)=\sum_{1\le r_1\le r_2\le\ldots\le r_l\le m}t_{j+1-r_1}t_{j+1-r_2}\ldots
t_{j+1-r_l}.
$$ 
Then, the moments of the B-spline $M_{j}$ of degree $m-2$,
with supp$(M_j)=[t_{j-m+1},t_{j}]$ and normalized by $\int M_j=1$, are given by (see e.g. \cite{nm}):
$$
\mu_j^{(l)}:=\int x^l  M_{j}(x)dx=(C_{m+l-1}^l)^{-1}\bar\sigma_l(T_j).
$$

\section{Uniformly bounded discrete quasi-interpolants exact on $\PP_2$}

It is possible to derive {\sl discrete} QIs from the de Boor-Fix QIs \cite{dbf} by replacing the values of derivatives $D^lf(\theta_i)/l!$ by divided differences at the points $\theta_i$ lying in
$\Sigma_i$. Doing this, we loose the property of projection on $S_m(I,T)$. However, by choosing
conveniently the divided differences, we can obtain some families of dQIs which are uniformly
bounded and exact on specific subspaces of polynomials.

Let us construct for example a family of dQIs of degree $m$ which are {\sl exact on} $\PP_2$.
We start from the de Boor-Fix  functionals truncated at order $2$:
$$
\lambda_j(f)=\frac{1}{m!}\sum_{l=0}^{2} (-1)^{m-l}D^{m-l}\psi_j(\tau)D^lf(\tau).
$$
where $\psi_j(t)=(t-t_{j-m+1})\ldots (t-t_j)$.
We obtain successively 
$$
D^{m}\psi_j(\tau)=(-1)^{m} m!,\;\;
D^{m-1}\psi_j(\tau)=(-1)^{m} m!(\tau-\theta_j),
$$
$$
D^{m-2}\psi_j(\tau)=\frac12 (-1)^{m} m!(\tau^2-2\theta_j t+\theta_j^{(2)}).
$$

More specifically, taking $\tau=\theta_j$, we get
$$
D^{m-1}\psi_j(\theta_j)=0,\;\; D^{m-2}\psi_j(\theta_j)=\frac12 (-1)^m
m!(\theta_j^{(2)}-\theta_j^2)
$$
Thus, we can we define the QI exact on $\PP_2$
$$
Q_2f=\sum_{j\in J} \lambda_j (f) B_j,
$$
whose coefficient functionals are given by
$$
\lambda_j (f)=f(\theta_j)-\frac12 \bar\theta_j^{(2)} D^2f(\theta_j),\;\;
{\rm with} \;\; \bar\theta_j^{(2)}=\theta_j^2-\theta_j^{(2)}.
$$
We recall the expansion (se e.g. \cite{db1}\cite{ms}):
$$
\bar\theta_j^{(2)}=\frac{1}{m^2 (m-1)}\sum_{1\le r<s\le m}(t_{j-r}-t_{j-s })^2.
$$
On the other hand, $\frac12 D^2f(\theta_j)$ coincides on the space $\PP_2$ with the
second order divided difference $[\theta_{j-1}, \theta_j,\theta_{j+1}]f$, 
therefore the dQI defined by
$$
Q_2^*f=\sum_{j\in J} \lambda_j^*(f) B_j,
$$
with coefficient functionals
$$
\lambda_j^*(f) =f(\theta_j)-\bar\theta_j^{(2)} [\theta_{j-1}, \theta_j,
\theta_{j+1}]f,
$$
is also exact on $\PP_2$. Moreover, one can write
$$
\lambda_j^*(f) =a_jf(\theta_{j-1})+b_j f(\theta_{j})+c_j f(\theta_{j+1}),\;\;{\rm with}
$$
$
a_j=-\bar\theta_j^{(2)}/\Delta \theta_{j-1}(\Delta \theta_{j-1}+\Delta \theta_{j}),\;\;
b_j=1+\bar\theta_j^{(2)}/\Delta \theta_{j-1}\Delta \theta_{j},\;\;
$
$c_j=-\bar\theta_j^{(2)}/\Delta \theta_{j}(\Delta \theta_{j-1}+\Delta \theta_{j}).$
So, according to the introduction
$$
\Vert Q_2^* \Vert_{\infty}\le \max_{j\in J} (\vert a_j \vert+\vert b_j \vert +\vert c_j \vert)
\le 1+2\max_{j\in J} \bar\theta_j^{(2)}/\Delta \theta_{j-1}\Delta \theta_{j}.
$$
The following theorem extends a result given for quadratic splines in \cite{mj}\cite{sab2}\cite{sab3}.

\begin{theorem}
For any degree $m$, the dQIs $Q_2^*$ are uniformly bounded. More specifically, for all
partitions of $I$:
$$
\Vert Q_2^* \Vert_{\infty}\le [\frac12 (m+4)].
$$
\end{theorem}
{\sl proof:} We only give the proof for $m=2k+1$, the case $m=2k$ being similar. For the sake of simplicity, we
take $j=m$, i.e. we shall determine an upper bound of the ratio
$$
N_m/D_m=\bar\theta_m^{(2)}/\Delta \theta_{m-1}\Delta \theta_{m}
$$
with
$$
N_m=\bar\theta_m^{(2)}=\frac{1}{m^2 (m-1)}\sum_{1\le r<s\le m} (t_r-t_s)^2.
$$
Setting $H=\sum_{i=2}^{m}h_i$, then we get a lower bound for the denominator
$$
D_m=\frac{1}{m^2}(t_{m}-t_0)(t_{m+1}-t_1)=\frac{1}{m^2}(h_1+H)(H+h_{m+1})\ge
\frac{H^2}{m^2}.
$$
The numerator $N_m$ is composed of $k$ pairs of sums $(S_p,S'_p)$
$$
S_p=\sum_{s-r=p} (t_r-t_s)^2,\;\; S'_p=\sum_{s-r=m-p} (t_r-t_s)^2,
$$
for  $1\le p\le k$.
Both sums contain at most $p$ times the terms $h_i^2$ and $2h_i h_j$ ($i\ne j$), hence we can
write
$S_p+S'_p\le 2p H^2$, which implies 
$$
N_m\le \frac{2H^2}{m^2(m-1)}(1+2+\ldots+k)=\frac{(k+1)H^2}{2m^2},
$$
so, we get
$$
N_m/D_m\le \frac12(k+1),
$$
and finally, for $m=2k+1$ odd
$$
\Vert Q_2^* \Vert_{\infty}\le k+2=\frac12 (m+3)=[\frac12 (m+4)].
$$
For $m=2k$, we obtain respectively $D_m\ge\frac{H^2}{m^2}$ and 
$N_m\le\frac{H^2}{4(m-1)}$, whence $N_m/D_m\le \frac{k^2}{2k-1}$, and finally
for $m=2k$ even
$$
\Vert Q_2^* \Vert_{\infty}\le k+2=\frac12 (m+4)=[\frac12 (m+4)].\quad \square
$$


\section{Existence and characterization of near-best discrete quasi-interpolants}

\subsection{ Existence of near-best dQIs}

We  consider the following family of dQIs defined on $I=\RR$ endowed with an
arbitrary non-uniform strictly increasing sequence of knots $T=\{t_i, \,i\in \ZZ\}$,
$$
Qf=Q_{p,q}f=\sum_{i\in\ZZ} \lambda_i(f) B_i.
$$ 
Their coefficient functionals depend on $2p+1$ parameters, with  $2p\ge m$,
$$
\lambda_i(f) =\sum_{s=-p}^{p} a_i(s)f(\theta_{i+s}),
$$
and they are exact on the space $\PP_q$, where $q\le m$.
The latter condition is equivalent to $Q e_r=e_r$ for all monomials of degrees $0\le r\le q$.
It implies that for all indices $i$, the parameters $a_i(s)$ satisfy the system of $q+1$
linear equations:
$$
\sum_{s=-p}^{p}a_i(s)\theta_{i+s}^r=\theta_i^{(r)},\quad 0\le r\le q.
$$
The Vandermonde matrix $V_i\in \RR^{(q+1)\times (2p+1)}$ of this system, with coefficients 
$V_i(r,s)=\theta_{i+s}^r$, is of maximal rank $q+1$, therefore there are $2p-q$ {\sl free parameters}. 
Denoting $b_i\in \RR^{q+1}$ the vector with components 
$b_i(r)=\theta_i^{(r)},\quad 0\le r\le q$, and by $a_i\in\RR^{2p+1}$ the vector with components 
$a_i(s)$, we consider the sequence of minimization problems, for $i\in \ZZ$:
\begin{equation}
\min \Vert a_i \Vert_1,\quad V_i a_i=b_i. \tag{$M_i$ }
\end{equation}
We have already seen in the introduction that $\nu_1^*(Q)=\max_{i\in\ZZ}\min\Vert a_i \Vert_1$ is an upper bound 
of $\Vert Q \Vert_{\infty}$ which is easier to evaluate than the true norm of the dQI.\\

\begin{theorem} The above minimization problems $(M_i)$ have always solutions, which, in general,
are non unique.\
\end{theorem}

{\sl proof:} 
The {\sl objective function being convex} and {\sl the domains being affine subspaces},
these classical optimization problems have always solutions, in general non
unique. $\square.$
\subsection{Characterization of optimal solutions}

For $b\in \RR^m$ and $A\in\RR^{m\times n}$, let us consider the $l_1$-minimization problem 
\begin{equation}
 \quad \min \Vert r(a)\Vert_1, \quad r(a)=b-Aa.\tag{M}
\end{equation}

We recall the characterization of optimal solutions for $l_1$-problems given in \cite{w}, chapter 6.
Define the sets
$$
Z(a)=\{1\le i\le m\,;\,r_i(a)=0\},
$$
$$
V(a)=\{v\in \RR^m \,;\,\Vert v\Vert_{\infty}\le 1, \; v_i={\rm sgn}(r_i(a)) \; {\rm for} \; i\notin Z(a) \}.
$$
\begin{theorem} The vector $a^*\in \RR^n$ is a solution of (M) if and only if there exists a vector $v^*\in V(a^*)$ 
satisfying $A^T v^*=0$.
\end{theorem}

\section{A  family of spline discrete quasi-interpolants exact on $\PP_2$}

In this section, we restrict our study to the subfamily of spline dQIs which
are exact on $\mathbb{P}_{2}$, i.e. we consider the dQIs $Q_{p}=Q_{p,2}$.

\noindent We shall need some set of indices. Let $\overline{K}=\left\{
-p,\ldots,p\right\} $, and $K^{\ast}=\left\{ -p,0,p\right\} $. Then, we can
write $K:=\overline{K}\setminus K^{\ast}=K_{1}\cup K_{2}$, where $%
K_{1}=\left\{ -p+1,\ldots,-1\right\} $, and $K_{2}=\left\{ 1,\ldots
,p-1\right\} $.

\noindent The three equations expressing the exactness of $Q_{p}$ on $%
\mathbb{P}_{2}$ can be written%
\begin{align*}
a_{i}\left( -p\right) +a_{i}\left( 0\right) +a_{i}\left( p\right) &
=1-\sum_{r\in K}a_{i}\left( r\right) \\
\theta_{i-p}a_{i}\left( -p\right) +\theta_{i}a_{i}\left( 0\right)
+\theta_{i+p}a_{i}\left( p\right) & =\theta_{i}-\sum_{r\in K}\theta
_{i+r}a_{i}\left( r\right) \\
\theta_{i-p}^{2}a_{i}\left( -p\right) +\theta_{i}^{2}a_{i}\left( 0\right)
+\theta_{i+p}^{2}a_{i}\left( p\right) & =\theta_{i}^{\left( 2\right)
}-\sum_{r\in K}\theta_{i+r}^{2}a_{i}\left( r\right)
\end{align*}

\noindent Let $\left( a_{i}^{\ast }\left( -p\right) ,a_{i}^{\ast }\left(
0\right) ,a_{i}^{\ast }\left( p\right) \right)$ be the unique solution of the
system with the right-hand side obtained by taking $a_{i}\left( r\right) =0$
for all $r\in K$. Using Cramer's rule, we obtain%
\begin{align*}
a_{i}^{\ast }\left( -p\right) & =-\overline{\theta }_{i}^{\left( 2\right)
}/\left( \theta _{i+p}-\theta _{i-p}\right) \left( \theta _{i}-\theta
_{i-p}\right) , \\
a_{i}^{\ast }\left( 0\right) & =1+\overline{\theta }_{i}^{\left( 2\right)
}/\left( \theta _{i+p}-\theta _{i}\right) \left( \theta _{i}-\theta
_{i-p}\right) , \\
a_{i}^{\ast }\left( p\right) & =-\overline{\theta }_{i}^{\left( 2\right)
}/\left( \theta _{i+p}-\theta _{i-p}\right) \left( \theta _{i+p}-\theta
_{i}\right) .
\end{align*}%
Then we can express the general solution of the above system in the form%
\begin{align*}
a_{i}\left( -p\right) & =a_{i}^{\ast }\left( -p\right) -\sum_{r\in
K_{1}}\alpha _{i,r}a_{i}\left( r\right) +\sum_{s\in K_{2}}\alpha
_{i,s}a_{i}\left( s\right) , \\
a_{i}\left( 0\right) & =a_{i}^{\ast }\left( 0\right) -\sum_{r\in K_{1}}\beta
_{i,r}a_{i}\left( r\right) -\sum_{s\in K_{2}}\beta _{i,s}a_{i}\left(
s\right) , \\
a_{i}\left( p\right) & =a_{i}^{\ast }\left( p\right) +\sum_{r\in
K_{1}}\gamma _{i,r}a_{i}\left( r\right) -\sum_{s\in K_{2}}\gamma
_{i,s}a_{i}\left( s\right) ,
\end{align*}%
with%
\begin{align*}
\alpha _{i,j}& =\left\{
\begin{array}{ll}
\dfrac{\left( \theta _{i}-\theta _{i+j}\right) \left( \theta _{i+p}-\theta
_{i+j}\right) }{\left( \theta _{i}-\theta _{i-p}\right) \left( \theta
_{i+p}-\theta _{i-p}\right) }, & \text{if }j\in K_{1}, \\
\dfrac{\left( \theta _{i+j}-\theta _{i}\right) \left( \theta _{i+p}-\theta
_{i+j}\right) }{\left( \theta _{i}-\theta _{i-p}\right) \left( \theta
_{i+p}-\theta _{i-p}\right) }, & \text{if }j\in K_{2},%
\end{array}%
\right. \  \\
\gamma _{i,j}& =\left\{
\begin{array}{ll}
\dfrac{\left( \theta _{i+j}-\theta _{i-p}\right) \left( \theta _{i}-\theta
_{i+j}\right) }{\left( \theta _{i+p}-\theta _{i-p}\right) \left( \theta
_{i+p}-\theta _{i}\right) }, & \text{if }j\in K_{1}, \\
\dfrac{\left( \theta _{i+j}-\theta _{i-p}\right) \left( \theta _{i+j}-\theta
_{i}\right) }{\left( \theta _{i+p}-\theta _{i-p}\right) \left( \theta
_{i+p}-\theta _{i}\right) }, & \text{if }j\in K_{2},%
\end{array}%
\right.
\end{align*}%
and%
\begin{equation*}
\beta _{i,j}=\frac{\left( \theta _{i+j}-\theta _{i-p}\right) \left( \theta
_{i+p}-\theta _{i+j}\right) }{\left( \theta _{i}-\theta _{i-p}\right) \left(
\theta _{i+p}-\theta _{i}\right) },\ j\in K_{1}\cup K_{2}.
\end{equation*}

\noindent We denote by $Q_{p}^{\ast}$ the spline dQI whose coefficient
functionals are%
\begin{equation*}
\lambda_{i}^{\ast}\left( f\right) =a_{i}^{\ast}\left( -p\right) f\left(
\theta_{i-p}\right) +a_{i}^{\ast}\left( 0\right) f\left( \theta _{i}\right)
+a_{i}^{\ast}\left( p\right) f\left( \theta_{i+p}\right) .
\end{equation*}
In that case, an upper bound of the norm of this QI is $\max_{i\in\mathbb{Z}%
}\nu_{i}^{\ast}$, where%
\begin{equation*}
\nu_{i}^{\ast}=\left\vert a_{i}^{\ast}\left( -p\right) \right\vert
+\left\vert a_{i}^{\ast}\left( 0\right) \right\vert +\left\vert a_{i}^{\ast
}\left( p\right) \right\vert .
\end{equation*}

\begin{theorem}
For all $p\geq m$, the infinity norms of the spline dQIs $Q_{p}^{\ast}$ are
uniformly bounded by $\frac{m+1}{m-1}$. This bound is independent of $p$ and
of the sequence of knots $T$.
\end{theorem}

\noindent\textit{proof}: We have to find a good upper bound of $%
\nu_{i}^{\ast }$. We recall that $\theta_{i}=\frac{1}{m}\sum_{r\in\mathbb{N}%
_{m}}t_{i+1-r}$ and $\overline{\theta}_{i}^{\left( 2\right) }=\frac{1}{%
m^{2}\left( m-1\right) }S_{1}$, with
\begin{equation*}
S_{1}=\sum_{1\leq r<s\leq m}\left( t_{i+1-r}-t_{i+1-s}\right) ^{2}.
\end{equation*}
Define the following sums:%
\begin{align*}
S_{2} & =\sum_{r\in\mathbb{N}_{m}}\left( t_{i+1-r}-t_{i+1-m-r}\right)
=\sum_{r\in\mathbb{N}_{m}}\sum_{k=1}^{m}h_{i+2-r-k}, \\
S_{3} & =\sum_{r\in\mathbb{N}_{m}}\left( t_{i+1+m-r}-t_{i+1-r}\right)
=\sum_{r\in\mathbb{N}_{m}}\sum_{k=1}^{m}h_{i+1-r+k}.
\end{align*}
As $p\geq m$, we obtain%
\begin{align*}
\theta_{i}-\theta_{i-p} & =\frac{1}{m}\sum_{r\in\mathbb{N}_{m}}\left(
t_{i+1-r}-t_{i+1-p-r}\right) =\frac{1}{m}\sum_{r\in\mathbb{N}_{m}}\sum
_{k=1}^{p}h_{i+2-r-k}\geq\frac{1}{m}S_{2}, \\
\theta_{i+p}-\theta_{i} & =\frac{1}{m}\sum_{r\in\mathbb{N}_{m}}\left(
t_{i+1+p-r}-t_{i+1-r}\right) =\frac{1}{m}\sum_{r\in\mathbb{N}_{m}}\sum
_{k=1}^{p}h_{i+1+k-r}\geq\frac{1}{m}S_{3}.
\end{align*}
The proof being essentially the same for all $i\in\mathbb{Z}$, we can
restrict our study to the case $i=m$. In that case, we get%
\begin{align*}
S_{2} & =mh_{1}+\sum_{k=1}^{m-1}k\left( h_{m+1-k}+h_{k+1-m}\right) \geq
S_{2}^{\prime}=\left( m-1\right) h_{2}+\cdots+2h_{m-1}+h_{m}, \\
S_{3} & =mh_{m+1}+\sum_{k=1}^{m-1}k\left( h_{2m+1-k}+h_{k+1}\right) \geq
S_{3}^{\prime}=\left( m-1\right) h_{m}+\cdots+2h_{3}+h_{2}.
\end{align*}
Setting $H_{k}=h_{2}+\cdots+h_{k+1}$, for $1\leq k\leq m-1$, and $H=H_{m-1}$
for short as in the proof of theorem 1, we have%
\begin{equation*}
S_{2}^{\prime}=\sum_{i=0}^{m-1}H_{i},\quad
S_{3}^{\prime}=H+\sum_{i=1}^{m-2}\left( H-H_{i}\right) =mH-S_{2}^{\prime},
\end{equation*}
whence%
\begin{equation*}
S_{2}S_{3}\geq S_{2}^{\prime}S_{3}^{\prime}=mH\sum_{i=0}^{m-1}H_{i}-\left(
\sum_{i=0}^{m-1}H_{i}\right) ^{2}.
\end{equation*}
Now, we come back to $S_{1}$ and we shall prove that $S_{1}\leq
S_{2}^{\prime }S_{3}^{\prime}\leq S_{2}S_{3}$. $S_{1}$ can be written under
the form%
\begin{equation*}
S_{1}=\sum_{r=1}^{m-1}\sum_{j=1}^{m-r}\left( h_{r+1}+\cdots+h_{r+j}\right)
^{2}=\sum_{i=1}^{m-1}H_{i}^{2}+\sum_{j=1}^{m-2}\sum_{i=j+1}^{m-1}\left(
H_{i}-H_{j}\right) ^{2},
\end{equation*}
from which we deduce%
\begin{equation*}
S_{1}=\left( m-1\right)
\sum_{i=1}^{m-1}H_{i}^{2}-2\sum_{j=1}^{m-2}H_{j}\sum_{i=j+1}^{m-1}H_{i}=m%
\sum_{i=1}^{m-1}H_{i}^{2}-\left( \sum_{i=1}^{m-1}H_{i}\right) ^{2}.
\end{equation*}
Then we use the fact that, for all $1\leq i\leq m-1$, $H_{i}\leq H$, whence $%
H_{i}^{2}\leq H\ H_{i}$ and $\sum_{i=0}^{m-1}H_{i}^{2}\leq
H\sum_{i=0}^{m-1}H_{i}$. So, we obtain%
\begin{equation*}
S_{1}\leq mH\sum_{i=0}^{m-1}H_{i}-\left( \sum_{i=0}^{m-1}H_{i}\right)
^{2}=S_{2}^{\prime}S_{3}^{\prime}\leq S_{2}S_{3}.
\end{equation*}
Finally, for all $i\in\mathbb{Z}$, we have%
\begin{equation*}
\nu_{i}^{\ast}=1+\frac{2}{m-1}\frac{S_{1}}{S_{2}S_{3}}\leq1+\frac{2}{m-1}=%
\frac{m+1}{m-1}\text{,}
\end{equation*}
whence $\displaystyle \left\Vert Q_{p}^{\ast}\right\Vert _{\infty}\leq\max_{i\in\mathbb{Z}%
}\nu_{i}^{\ast}\leq\frac{m+1}{m-1}$. $\square$\\

In the next section we prove that the quasi-interpolants $Q_{p}^{\ast}$ are
near-best in the sense of section 4 under some additional conditions on the partitions.


\section{The family $Q_{p}^*$  of   discrete quasi-interpolants is near-best}

Let us write the minimization problem ($P_d$) of section 4 in
Watson's form. Taking into account the expression of the solution $a_{i}$ of
the system equivalent to the exactness on $\mathbb{P}_{2}$ of $Q_{p,2}$, we
can write%
\begin{equation*}
\left\Vert a_{i}\right\Vert _{1}=\left\Vert a_{i}^{\ast }-A_{i}\widetilde{a}%
_{i}\right\Vert ,
\end{equation*}%
where%
\begin{align*}
\widetilde{a}_{i}& =\left( a_{i}\left( -p+1\right) ,\ldots ,a_{i}\left(
-1\right) ,a_{i}\left( 1\right) ,\ldots ,a_{i}\left( p-1\right) \right)
^{T}\in \mathbb{R}^{2p-2}, \\
a_{i}^{\ast }& =\left( a_{i}^{\ast }\left( -p\right) ,0,\ldots
,0,a_{i}^{\ast }\left( 0\right) ,0,\ldots ,a_{i}^{\ast }\left( p\right)
\right) ^{T}\in \mathbb{R}^{2p+1},
\end{align*}%
and $A_{i}\in \mathbb{R}^{\left( 2p+1\right) \times \left( 2p-2\right) }$ is
given by%
\begin{equation*}
A_{i}=\left(
\begin{array}{cccccc}
\alpha _{i,-p+1} & \cdots  & \alpha _{i,-1} & -\alpha _{i,-1} & \cdots  &
\alpha _{i,p-1} \\
-1 & \cdots  & 0 & 0 & \cdots  & 0 \\
\vdots  & \ddots  & \vdots  & \vdots  & \ddots  & \vdots  \\
0 & \cdots  & -1 & 0 & \cdots  & 0 \\
\beta _{i,-p+1} & \cdots  & \beta _{i,-1} & \beta _{i,1} & \cdots  & \beta
_{i,p-1} \\
0 & \cdots  & 0 & -1 & \cdots  & 0 \\
\vdots  & \ddots  & \vdots  & \vdots  & \ddots  & \vdots  \\
0 & \cdots  & 0 & 0 & \cdots  & -1 \\
-\gamma _{i,-p+1} & \cdots  & -\gamma _{i,-1} & \gamma _{i,1} & \cdots  &
\gamma _{i,p-1}%
\end{array}%
\right) .
\end{equation*}

\begin{theorem}
Assume that the sequence of knots $T$ satisfies, for all $i\in \mathbb{Z}$,
the following properties%
\begin{equation*}
\theta _{i-1}+\theta _{i}\leq \theta _{i-p}+\theta _{p}\leq \theta
_{i}+\theta _{i+1,}
\end{equation*}%
then, for all $i\in \mathbb{Z}$, $a_{i}^{\ast }$ is an optimal solution of
the local minimization problem $(M_i)$. Thus, for all $p\geq m$, the
spline dQIs $Q_{p}^{\ast }$ (theorem 4) are near-best.
\end{theorem}

\noindent\noindent\textit{proof}: According to theorem 3,
 we must find a vector $v^{\ast}\in\mathbb{R}^{2p+1}$ satisfying%
\begin{equation*}
\left\Vert v^{\ast}\right\Vert _{\infty}\leq1,\quad
A_{i}^{T}v^{\ast}=0,\quad v^{\ast}\left( r\right) ={\rm sgn}\left(
a_{i}^{\ast}\left( r\right) \right) \text{ for }r=-p,0,p\text{.}
\end{equation*}
Let us choose $v^{\ast}\left( -p\right) =-1$, $v^{\ast}\left( 0\right) =1$, $%
v^{\ast}\left( p\right) =-1$, and%
\begin{equation*}
v^{\ast}\left( j\right) =\left\{
\begin{array}{ll}
-\alpha_{i,r}+\beta_{i,r}+\gamma_{i,r}, & \text{if }j\in K_{1}, \\
\alpha_{i,r}+\beta_{i,r}-\gamma_{i,r}, & \text{if }j\in K_{2}.%
\end{array}
\right.
\end{equation*}
Then it is easy to verify that the equations $A_{i}^{T}v^{\ast}=0$ are
satisfied. Moreover, the above expressions of $a_{i}^{\ast}\left( r\right) $
for $r=-p,0,p$, with $\overline{\theta}_{i}^{\left( 2\right) }>0$ imply that
$v^{\ast}\left( r\right) ={\rm sgn}\left( a_{i}^{\ast}\left( r\right)
\right) $. It only remains to prove that $\left\vert v^{\ast }\left(
j\right) \right\vert \leq1$ for $j\in K_{1}\cup K_{2}$. As $%
\beta_{i,r}=1-\alpha_{i,r}+\gamma_{i,r}$ for $r\in K_{1}$ and $\beta
_{i,s}=1+\alpha_{i,s}-\gamma_{i,s}$ for $s\in K_{2}$, it is equivalent to
prove%
\begin{equation*}
0\leq\alpha_{i,r}-\gamma_{i,r}\leq1,\quad0\leq\gamma_{i,s}-\alpha_{i,s}%
\leq1,\quad\text{for }\left( r,s\right) \in K_{1}\times K_{2}.
\end{equation*}
We only detail the proof for $r\in K_{1}$, that for $s\in K_{2}$ being quite
similar. Using the expressions of $\alpha_{i,r}$ and $\gamma_{i,r}$ given in
section 5, we get%
\begin{equation*}
\alpha_{i,r}-\gamma_{i,r}=\frac{\left( \theta_{i}-\theta_{i+r}\right) \left[
\left( \theta_{i+p}+\theta_{i-p}\right) -\left(
\theta_{i+r}+\theta_{i}\right) \right] }{\left(
\theta_{i}-\theta_{i-p}\right) \left( \theta_{i+p}-\theta_{i}\right) },
\end{equation*}
and we shall have $\alpha_{i,r}-\gamma_{i,r}\geq0$ if and only if%
\begin{equation*}
\theta_{i+r}+\theta_{i}\leq\theta_{i+p}+\theta_{i-p}
\end{equation*}
for all $r\in K_{1}$. However, since we have $\theta_{i+r}+\theta_{i}\leq%
\theta_{i-1}+\theta_{i}$, there only remains the unique condition%
\begin{equation*}
\theta_{i-1}+\theta_{i}\leq\theta_{i-p}+\theta_{i+p}.
\end{equation*}
The other inequality $\alpha_{i,r}-\gamma_{i,r}\leq1$ can be written%
\begin{equation*}
\left( \theta_{i}-\theta_{i+r}\right) \left[ \left(
\theta_{i+p}+\theta_{i-p}\right) -\left( \theta_{i+r}+\theta_{i}\right) %
\right] \leq\left( \theta_{i}-\theta_{i-p}\right) \left( \theta_{i+p}-\theta
_{i}\right) .
\end{equation*}
Setting $\delta_{1}=\theta_{i+r}-\theta_{i-p}$, $\delta_{2}=\theta_{i}-%
\theta_{i+r}$, and $\delta_{3}=\theta_{i+p}-\theta_{i}$, the latter
inequality can be written $\delta_{2}\left( \delta_{3}-\delta_{1}\right)
\leq\delta_{3}\left( \delta_{2}+\delta_{1}\right) $, or equivalently $%
\delta_{1}\left( \delta_{2}+\delta_{3}\right) \geq0$ which is obviously
satisfied. For $s\in K_{2}$, the inequalities $0\leq\gamma_{i,s}-\alpha
_{i,s}\leq1$ are satisfied if and only if%
\begin{equation*}
\theta_{i-p}+\theta_{i+p}\leq\theta_{i}+\theta_{i+1},
\end{equation*}
whence the conditions on the sequence of knots. $\square$

\begin{remark}
Theorem 5 imposes some additional conditions on the sequence of knots. For quadratic splines,
we have studied arithmetic and geometric sequences: in both cases, the higher is $p$, 
the stronger are the conditions and, for $p\to +\infty$, $T$ is closer and closer to a uniform sequence.
\end{remark}

\begin{remark}
Even if the partition $T$ does not satisfy the hypotheses of theorem 5, the operator $Q^*_p$ is still a good QI because its infinity norm is small and uniformly bounded.

\end{remark}

-
\section{Uniformly bounded integral quasi-interpolants of Goodman-Sharma type}
General integral spline quasi-interpolants (iQIs) already appear in \cite{db1}-\cite{db2}\cite{ls}\cite{sab1}\cite{ss}. Here we have chosen to study a family of QIs which we call Goodman-Sharma (GS-) type
iQIs, as they first appear in  \cite{gs}. They seem simpler and more interesting
than those studied in  \cite{sab1} and  \cite{ss}. In the case of splines of
degree $m$ on a partition in $n$ subintervals of a bounded interval $I=\left[
a,b\right] $, the simplest GS-type iQI can be written as follows:%
\begin{equation*}
G_{1}f=f\left( t_{0}\right) B_{0}+\sum_{i=1}^{n+m-2}\lambda_{i}\left(
f\right) B_{i}+f\left( t_{n}\right) B_{n+m-1},
\end{equation*}
where the integral coefficient functionals are defined by%
\begin{equation*}
\lambda_{i}\left( f\right) =\int_{a}^{b}M_{i}\left( t\right) f\left(
t\right) dt=\left\langle M_{i},f\right\rangle ,
\end{equation*}
$M_{i}$ being the B-spline of degree $m-2$ with support $\Sigma_{i}=\left[
t_{i-m+1},t_{i}\right] $, normalized by $\lambda_{i}\left( e_{0}\right)
=\mu_{i}^{\left( 0\right) }=\int_{\mathbb{R}}M_{i}=1$. It is easy to verify
that $G_{1}$ is exact on $P_{1}$ and that $\left\Vert G_{1}\right\Vert
_{\infty}=1$. In this section, we shall study the family of GS-type iQIs
defined by%
\begin{equation*}
G_{2}f=f\left( t_{0}\right) B_{0}+\sum_{i=1}^{n+m-2}\left[
a_{i}\lambda_{i-1}\left( f\right) +b_{i}\lambda_{i}\left( f\right)
+c_{i}\lambda_{i+1}\left( f\right) \right] B_{i}+f\left( t_{n}\right)
B_{n+m-1},
\end{equation*}
which we want to be exact on $\mathbb{P}_{2}$. The three constraints $%
G_{2}e_{k}=e_{k}$, $k=0,1,2$, lead to the following system of equations, for
each $1\leq i\leq n+m-2$:%
\begin{equation*}
a_{i}+b_{i}+c_{i}=1,\
\theta_{i-1}a_{i}+\theta_{i}b_{i}+\theta_{i+1}c_{i}=\theta_{i},\
\mu_{i-1}^{\left( 2\right) }a_{i}+\mu_{i}^{\left( 2\right)
}b_{i}+\mu_{i+1}^{\left( 2\right) }c_{i}=\theta_{i}^{\left( 2\right) }.
\end{equation*}
We recall the values of the first moments of $M_{i}$:%
\begin{align*}
\mu_{i}^{\left( 1\right) } & =\frac{1}{m}\sum_{1\leq r\leq m}t_{i+1-r}=\frac{%
1}{m}\sum_{1\leq r\leq m}t_{i-m+r}=\theta_{i}, \\
\mu_{i}^{\left( 2\right) } & =\frac{2}{m\left( m+1\right) }\sum_{1\leq r\leq
s\leq m}t_{i+1-r}t_{i+1-s}=\frac{2}{m\left( m+1\right) }\sum_{1\leq r\leq
s\leq m}t_{i-m+r}t_{i-m+s}.
\end{align*}

\begin{theorem}
For any degree $m$, the iQIs $G_{2}$ are uniformly bounded indepently of the
partition of $I$. For $m=2k$ or $2k+1$, there holds%
\begin{equation*}
\left\Vert G_{2}\right\Vert _{\infty}\leq2k+3.
\end{equation*}
\end{theorem}

\noindent \textit{proof}: Taking the differences of the second and third
equations above ($G_{2}e_{k}=e_{k}$, $k=1,2$) with the first one ($%
G_{2}e_{0}=e_{0}$) times $\theta _{i}$ and $\mu _{i}^{\left( 2\right) }$
resp., we get%
\begin{equation*}
\left( \theta _{i}-\theta _{i-1}\right) a_{i}=\left( \theta _{i+1}-\theta
_{i}\right) c_{i},\quad \left( \mu _{i}^{\left( 2\right) }-\mu
_{i-1}^{\left( 2\right) }\right) a_{i}+\left( \mu _{i}^{\left( 2\right)
}-\mu _{i+1}^{\left( 2\right) }\right) c_{i}=\mu _{i}^{\left( 2\right)
}-\theta _{i}^{\left( 2\right) }.
\end{equation*}%
Setting $a_{i}=\left( \theta _{i+1}-\theta _{i}\right) \alpha _{i}$ and $%
c_{i}=\left( \theta _{i}-\theta _{i-1}\right) \alpha _{i}$, we obtain%
\begin{equation}
\left[ \Delta \mu _{i-1}^{\left( 2\right) }\Delta \theta _{i}-\Delta \mu
_{i}^{\left( 2\right) }\Delta \theta _{i-1}\right] \alpha _{i}=\mu
_{i}^{\left( 2\right) }-\theta _{i}^{\left( 2\right) }.  \tag{E}
\end{equation}%
Using the expressions of the various coefficients in terms of symmetric
functions of the knots, we obtain the following form for the coefficient of $%
\alpha _{i}$ in equation~{(E)}%
\begin{equation*}
\frac{2}{m^{2}\left( m+1\right) }\left[ \left( t_{i+1}-t_{i+1-m}\right)
\Delta \overline{\sigma }_{2}\left( T_{i-1}\right) -\left(
t_{i}-t_{i-m}\right) \Delta \overline{\sigma }_{2}\left( T_{i}\right) \right]
.
\end{equation*}%
Setting $l_{i}=t_{i}-t_{i+1-m}$, we can write $t_{i}-t_{i-m}=h_{i+1-m}+l_{i}$%
, $t_{i+1}-t_{i-m}=h_{i+1-m}+l_{i}+h_{i+1}$ and $%
t_{i+1}-t_{i+1-m}=l_{i}+h_{i+1}$. Then we have%
\begin{equation*}
\Delta \overline{\sigma }_{2}\left( T_{i}\right) =\left(
t_{i+1}-t_{i+1-m}\right) \sum_{r=i-m+1}^{i+1}t_{r},
\end{equation*}%
and the coefficient of $\alpha _{i}$ in equation~{(E)} is given by%
\begin{align*}
\Delta \mu _{i-1}^{\left( 2\right) }\Delta \theta _{i}-\Delta \mu
_{i}^{\left( 2\right) }\Delta \theta _{i-1}& =-\left(
t_{i+1}-t_{i+1-m}\right) \left( t_{i}-t_{i-m}\right) \left(
t_{i+1}-t_{i-m}\right)  \\
& =-\left( l_{i}+h_{i+1}\right) \left( h_{i+1-m}+l_{i}\right) \left(
h_{i+1-m}+l_{i}+h_{i+1}\right) .
\end{align*}%
Now, writting $\overline{\sigma }_{2}\left( T_{i}\right) =\sigma _{2}\left(
T_{i}\right) +\sum_{r=i+1-m}^{i}t_{r}^{2}$, we get successively%
\begin{align*}
\mu _{i}^{\left( 2\right) }-\theta _{i}^{\left( 2\right) }& =\frac{2}{%
m\left( m+1\right) }\overline{\sigma }_{2}\left( T_{i}\right) -\frac{2}{%
m\left( m-1\right) }\sigma _{2}\left( T_{i}\right)  \\
& =\frac{2}{m\left( m^{2}-1\right) }\left( \left( m-1\right) \overline{%
\sigma }_{2}\left( T_{i}\right) -\left( m+1\right) \sigma _{2}\left(
T_{i}\right) \right)  \\
& =\frac{2}{m\left( m^{2}-1\right) }\left( \left( m-1\right)
\sum_{r=i+1-m}^{i}t_{r}^{2}-2\sigma _{2}\left( T_{i}\right) \right)  \\
& =\frac{2}{m\left( m^{2}-1\right) }\sum_{i+1-m\leq r<s\leq i}\left(
t_{r}-t_{s}\right) ^{2}.
\end{align*}%
Setting $\omega _{i}=\sum_{i+1-m\leq r<s\leq i}\left( t_{r}-t_{s}\right) ^{2}
$, we obtain%
\begin{equation*}
\alpha _{i}=-\frac{m}{m-1}\frac{\omega _{i}}{\left( h_{i+1-m}+l_{i}\right)
\left( h_{i+1-m}+l_{i}+h_{i+1}\right) \left( l_{i}+h_{i+1}\right) },
\end{equation*}%
from which we deduce%
\begin{align*}
a_{i}& =-\frac{1}{m-1}\frac{\omega _{i}}{\left( h_{i+1-m}+l_{i}\right)
\left( h_{i+1-m}+l_{i}+h_{i+1}\right) }, \\
c_{i}& =-\frac{1}{m-1}\frac{\omega _{i}}{\left(
h_{i+1-m}+l_{i}+h_{i+1}\right) \left( l_{i}+h_{i+1}\right) },
\end{align*}%
and%
\begin{equation*}
b_{i}=1-a_{i}-c_{i},\quad \left\vert a_{i}\right\vert +\left\vert
b_{i}\right\vert +\left\vert c_{i}\right\vert =1+2\left( \left\vert
a_{i}\right\vert +\left\vert c_{i}\right\vert \right) .
\end{equation*}%
On the other hand, we have%
\begin{equation*}
\left\vert a_{i}\right\vert +\left\vert c_{i}\right\vert =\frac{1}{m-1}\frac{%
\omega _{i}\left( h_{i+1-m}+2l_{i}+h_{i+1}\right) }{\left(
h_{i+1-m}+l_{i}\right) \left( h_{i+1-m}+l_{i}+h_{i+1}\right) \left(
l_{i}+h_{i+1}\right) }.
\end{equation*}%
As $\omega _{i}\leq k\left( k+1\right) l_{i}^{2}$ for $m=2k+1$ (see proof of
theorem~1 with respect the upper bound for $N_{m}$), we obtain%
\begin{equation*}
\omega _{i}\leq k\left( k+1\right) \left( h_{i+1-m}+l_{i}\right) \left(
l_{i}+h_{i+1}\right) .
\end{equation*}%
Moreover, it is obvious that%
\begin{equation*}
h_{i+1-m}+2l_{i}+h_{i+1}\leq 2\left( h_{i+1-m}+l_{i}+h_{i+1}\right) .
\end{equation*}%
Therefore, we finally obtain for all $i$%
\begin{equation*}
\left\vert a_{i}\right\vert +\left\vert c_{i}\right\vert \leq \frac{2k\left(
k+1\right) }{m-1}=k+1=\frac{1}{2}\left( m+1\right) ,
\end{equation*}%
whence the uniform upper bound for the norm of $G_{2}$%
\begin{equation*}
\left\Vert G_{2}\right\Vert _{\infty }\leq m+2,
\end{equation*}%
which is both independent of the partition and of the (odd) degree of the
spline. For $m=2k$ even, a similar computation leads to the uniform bound%
\begin{equation*}
\left\Vert G_{2}\right\Vert _{\infty }\leq m+3.
\end{equation*} $\square$


\section{Existence and characterization of near-best integral quasi-interpolants}
Now, we consider the family of iQIs%
\begin{equation*}
G_{p,q}f=\sum_{i\in\mathbb{Z}}\lambda_{i}\left( f\right) B_{i},
\end{equation*}
whose coefficient functionals depend on $2p+1$ parameters:%
\begin{equation*}
\lambda_{i}\left( f\right) =\sum_{s=-p}^{p}a_{i}\left( s\right)
\int_{\Sigma_{i+s}}M_{i+s}\left( t\right) f\left( t\right) dt,
\end{equation*}
and which are exact on $\mathbb{P}_{q}$. As in section 7, $M_{j}$ is the
B-spline of degree $m-2$, with support $\Sigma_{j}=\left[ t_{j+1-m},t_{j}%
\right] $ normalized by $\int_{\mathbb{R}}M_{j}=1$. The constraints $%
G_{p,q}\left( e_{r}\right) =e_{r}$ are equivalent to the following systems
of linear equations, for all $i\in\ZZ$:%
\begin{equation*}
\sum_{s=-p}^{p}\mu_{i+s}^{\left( r\right) }a_{i}\left( s\right)
=\theta_{i}^{\left( r\right) },\ 0\leq r\leq q,
\end{equation*}
whose matrix coefficients $W_{i}\in\mathbb{R}^{\left( q+1\right)
\times\left( 2p+1\right) }$, defined by $W_{i}\left( r,s\right) =\mu
_{i+s}^{\left( r\right) }$, is of maximal rank $q+1$ (see e.g.  \cite{ss}), thus
there remains $2p-q$ free parameters. In view of the introduction, as $%
\max_{i\in\mathbb{Z}}\left\Vert a_{i}\right\Vert _{1}$ is an upper bound of
the true infinity norm of the iQI, we want to solve the minimization
problems, for all $i$:%
\begin{equation}
\min\left\Vert a_{i}\right\Vert _{1},\quad W_{i}a_{i}=b_{i}.
\tag{$\overline M_i$}
\end{equation}
As in section 4.1, the objective function being convex and the domains being
affine subspaces, we can conclude

\begin{theorem}
The above minimization problems~$(\overline M_i)$ have always solutions, in
general non unique.
\end{theorem}

\noindent As in section 4.2, we shall use the characterizacion of optimal
solutions for $l_{1}$-problems (M) given in  \cite{w}, chapter 6.


\section{A  family of integral spline quasi-interpolants exact on $\PP_2$}

In this section, we restrict our study to the subfamily $q=2$ of the above
spline iQIs which are exact on $\mathbb{P}_{2}$. Moreover, we assume that $%
p\geq m$ in order to insure that the three sets of knots $T_{i-p}$, $T_{i}$,
and $T_{i+p}$ are pairwise disjoint. Now, the matrix coefficients of the
linear system equivalent to the exactness of $G_{p,2}$ on $\mathbb{P}_{2}$
is of maximal rank $3$, and we have $2p-2$ free parameters. Let us denote by
$G_{p}^{\ast}$ the spline iQI whose coefficient functionals are%
\begin{equation*}
\lambda_{i}^{\ast}\left( f\right) =a_{i}^{\ast}\left( -p\right) \left\langle
M_{i-p},f\right\rangle +a_{i}^{\ast}\left( 0\right) \left\langle
M_{i},f\right\rangle +a_{i}^{\ast}\left( p\right) \left\langle
M_{i+p},f\right\rangle ,
\end{equation*}
where%
\begin{align*}
a_{i}^{\ast}\left( -p\right) +a_{i}^{\ast}\left( 0\right) +a_{i}^{\ast
}\left( p\right) & =1, \\
\theta_{i-p}a_{i}^{\ast}\left( -p\right) +\theta_{i}a_{i}^{\ast}\left(
0\right) +\theta_{i+p}a_{i}^{\ast}\left( p\right) & =\theta_{i}, \\
\mu_{i-p}^{\left( 2\right) }a_{i}^{\ast}\left( -p\right) +\mu_{i}^{\left(
2\right) }a_{i}^{\ast}\left( 0\right) +\mu_{i+p}^{\left( 2\right)
}a_{i}^{\ast}\left( p\right) & =\theta_{i}^{\left( 2\right) },
\end{align*}
that is, their coefficients are the unique solution of the system obtained
by taking $a_{i}\left( r\right) =0$ for all $r\in K$ (we use the same
notations as in section~5).

\noindent The following result is the analog of theorem~4 for iQIs.

\begin{theorem}
For any degree $m\ge 2$, and for all $p\geq m$, the iQIs $G_{p}^{\ast}$ are
uniformly bounded independently of the partition of $I$. More specifically,
there holds 
\begin{equation*}
\left\Vert G_{p}^{\ast}\right\Vert _{\infty}\leq1+\frac{1}{4}C\left( m\right),\;\;
with\;\;
C\left( m\right) =\left\{
\begin{array}{ll}
\frac{m^{2}(m+2)}{\left( m-1\right) ^{2}} & \text{for }m\text{ even,} \\
\frac{\left( m+1\right)^2 }{m-1} & \text{for }m\text{ odd.}%
\end{array}
\right.
\end{equation*}
\end{theorem}

\noindent\textit{proof}: We shall prove that $\left\vert a_{i}^{\ast}\left(
-p\right) \right\vert +\left\vert a_{i}^{\ast}\left( 0\right) \right\vert
+\left\vert a_{i}^{\ast}\left( p\right) \right\vert \leq1+\frac{1}{4}C\left(
m\right) $ for all $i\in\mathbb{Z}$, which is sufficient to insure the
result. For the sake of simplicity, we can assume that $i=m$. Solving the
corresponding linear system, we get%
\begin{equation*}
a_{m}^{\ast}\left( -p\right) =-\frac{\xi_{m}\Delta_{p}\theta_{m}}{\delta _{m}%
},\quad a_{m}^{\ast}\left( p\right) =-\frac{\xi_{m}\Delta_{p}\theta_{m-p}}{%
\delta_{m}},\quad a_{0}^{\ast}\left( 0\right) =1-a_{m}^{\ast }\left(
-p\right) -a_{m}^{\ast}\left( p\right) ,
\end{equation*}
where%
\begin{equation*}
\xi_{m}=\mu_{m}^{\left( 2\right) }-\theta_{m}^{\left( 2\right) },\quad\text{%
and \quad}\delta_{m}=\Delta_{p}\theta_{m-p}\Delta_{p}\mu _{m}^{\left(
2\right) }-\Delta_{p}\theta_{m}\Delta_{p}\mu_{m-p}^{\left( 2\right) },
\end{equation*}
with%
\begin{equation*}
\Delta_{p}\theta_{l}=\theta_{l+p}-\theta_{l},\quad\Delta_{p}\mu_{l}^{\left(
2\right) }=\mu_{l+p}^{\left( 2\right) }-\mu_{l}^{\left( 2\right) }\,\,,\quad
l=m,m-p.
\end{equation*}
As proved in theorem 6, we have%
\begin{equation*}
\xi_{m}=\frac{2}{m\left( m^{2}-1\right) }\omega_{m}=\frac{2}{m\left(
m^{2}-1\right) }\sum_{1\leq r<s\leq m}\left( t_{r}-t_{s}\right) ^{2}.
\end{equation*}
For the expression of $\delta_{m}$, we need some additional notations. 
For $1\leq i\leq m$, we define 
\begin{equation*}
l_{i}=t_{1}-t_{i-p},\quad l_{i}^{\prime}=t_{i}-t_{i-p}=\sum_{r=2}^{i}h_{r}+l_{i}.
\end{equation*}
Similarly, for $1\leq j\leq m$, we define%
\begin{equation*}
l_{m+j}=t_{m+p}-t_{j},\quad l_{m+j}^{\prime}=t_{j+p}-t_{j}=l_{m+j}+\sum_{s=j}^{m-1}h_{s+1}.
\end{equation*}
Taking into account the definitions of $\theta_{i}$ and $\mu_{i}^{\left(
2\right) }$ for $i=m-p,m,m+p$, we obtain after some algebraic calculations
the following expressions:%
\begin{align*}
\Delta_{p}\theta_{m-p} & =\frac{1}{m}\sum_{i=1}^{m}l_{i}^{\prime},\quad%
\Delta_{p}\theta_{m}=\frac{1}{m}\sum_{j=1}^{m}l_{m+j}^{\prime}, \\
\Delta_{p}\mu_{m-p}^{\left( 2\right) } & =\frac{2}{m\left( m+1\right) }%
\sum_{i=1}^{m}l_{i}^{\prime}\left(
\sum_{r=i}^{m}t_{r-p}+\sum_{s=1}^{i}t_{s}\right) , \\
\Delta_{p}\mu_{m}^{\left( 2\right) } & =\frac{2}{m\left( m+1\right) }%
\sum_{j=1}^{m}l_{m+j}^{\prime}\left(
\sum_{r=j}^{m}t_{r}+\sum_{s=1}^{j}t_{p+s}\right) .
\end{align*}
Now we can write%
\begin{equation*}
a_{m}^{\ast}\left( -p\right) =-\frac{1}{m-1}\frac{\omega_{m}}{D_{m}}%
\sum_{j=1}^{m}l_{m+j}^{\prime},\quad a_{m}^{\ast}\left( p\right) =-\frac {1}{%
m-1}\frac{\omega_{m}}{D_{m}}\sum_{i=1}^{m}l_{i}^{\prime},
\end{equation*}
where $D_{m} $ is equal to%
$$
\left(\sum_{i=1}^{m}l_{i}^{\prime}\right)\left( \sum_{j=1}^{m}l_{m+j}^{\prime
}\left( \sum_{r=j}^{m}t_{r}+\sum_{s=1}^{j}t_{p+s}\right) \right)
-\left(\sum_{j=1}^{m}l_{m+j}^{\prime}\right)\left( \sum_{i=1}^{m}l_{i}^{\prime}\left(
\sum_{r=i}^{m}t_{r-p}+\sum_{s=1}^{i}t_{s}\right) \right)
$$
$$
=\sum_{i=1}^{m}\sum_{j=1}^{m}l_{i}^{\prime}l_{m+j}^{\prime}\left(
\sum_{r=j}^{m}t_{r}+\sum_{s=1}^{j}t_{p+s}-\sum_{r=i}^{m}t_{r-p}-%
\sum_{s=1}^{i}t_{s}\right) .
$$
Therefore,  as we shall see later (p.18) than $D_{m}>0$, we obtain 
\begin{equation*}
\left\vert a_{m}^{\ast}\left( -p\right) \right\vert +\left\vert a_{m}^{\ast
}\left( p\right) \right\vert =\frac{1}{m-1}\frac{\omega_{m}}{D_{m}}\left(
\sum_{j=1}^{m}l_{m+j}^{\prime}+\sum_{i=1}^{m}l_{i}^{\prime}\right) =\frac {1%
}{m-1}\frac{N_{m}}{D_{m}},
\end{equation*}
where%
\begin{equation*}
N_{m}=\left(\sum_{k=1}^{2m}l_{k}^{\prime}\right)\omega_{m}=\left(\sum_{k=1}^{2m}l_{k}^{\prime
}\right)\sum_{1\leq r<s\leq m}\left( t_{r}-t_{s}\right) ^{2}.
\end{equation*}
Let us compute an upper bound for $N_m$.
Let $H=\sum_{i=2}^{m}h_{i}$, as in the proof of  theorem 1. In the
first sum,  we have $l_{i}^{\prime}+l_{m+i}^{\prime
}=l_{i}+l_{i+m}+H$, for $1\leq i\leq m$, so we can write%
\begin{equation*}
\sum_{k=1}^{2m}l_{k}^{\prime}=\sum_{k=1}^{2m}l_{k}+mH.
\end{equation*}
For $\omega_{m}$, we have already seen (also in the proof of  theorem 1) that%
\begin{equation*}
\omega_{m}\le c\left( m\right) H^{2},
\end{equation*}
where $c\left( m\right) =k^{2}$ for $m=2k$ and $c\left( m\right) =k\left(
k+1\right) $ for $m=2k+1$. So, we finally  obtain the following  upper bound for $N_{m}$:%
\begin{equation*}
N_{m}\leq c\left( m\right) H^{2}\left( \sum_{k=1}^{2m}l_{k}+mH\right) .
\end{equation*}
Now, we will compute a lower bound for $D_{m}$. Let
\begin{equation*}
L_{i,m+j}=\sum_{r=j}^{m}t_{r}+%
\sum_{s=1}^{j}t_{p+s}-\sum_{r=i}^{m}t_{r-p}-\sum_{s=1}^{i}t_{s}
\end{equation*}
be the coefficient of $l_{i}^{\prime}l_{m+j}^{\prime}$ in the double sum defining $D_m$.
Therefore we have
$$
D_m\ge l'_ml'_{m+1}L_{m,m+1}+\sum_{i=1}^{m-1} l'_il'_{m+1}L_{i,m+1}+\sum_{j=2}^ml'_ml'_{m+j}L_{m,m+j}.
$$
We first observe that%
\begin{equation*}
l'_ml'_{m+1}L_{m,m+1}=\left( l_{m}+H\right) \left( l_{m+1}+H\right) \left(
l_{m}+l_{m+1}+H\right) \geq2H^{2}\left( l_{m}+l_{m+1}\right) +H^{3}.
\end{equation*}
Then, we obtain successively,  for $1\leq i\leq m-1$,%
\begin{equation*}
L_{i,m+1}=
\sum_{r=i+1}^{m}t_{r}+t_{p+1}-\sum_{r=i}^{m}t_{r-p}\ge \left( t_{m}-t_{m-1-p}\right) +\left(
t_{p+1}-t_{m-p}\right) \ge 2H.
\end{equation*}
Similarly, for $2\leq j\leq m$,%
\begin{equation*}
L_{m,m+j}=
\sum_{s=1}^{j}t_{p+s}-t_{m-p}-\sum_{r=1}^{j-1}t_{r}\ge \left( t_{p+1}-t_{m-p}\right) +\left(
t_{p+2}-t_{1}\right) \geq2H.
\end{equation*}
On the other hand, we also have successively
$$
l'_il'_{m+1}=(l_i+\sum_{r=2}^i h_r)(l_{m+1}+H)\ge H(l_i+\sum_{r=2}^i h_r),
$$
$$
l'_{m}l'_{m+j}=(l_m+H)(l_{m+j}+\sum_{s=j+1}^m h_s)\ge H(l_{m+j}+\sum_{s=j+1}^m h_s).
$$
From these inequalities, we deduce%
\begin{equation*}
D_{m}\geq H^3+2H^2\left(\sum_{k=1}^{2m}l_k+\sum_{i=2}^{m-1}\sum_{r=2}^i h_r+\sum_{j=2}^{m-1}
\sum_{s=j+1}^{m-1}h_s \right).
\end{equation*}
Now, it is easy to see that 
$$
\sum_{i=2}^{m-1}\sum_{r=2}^i h_r+\sum_{j=2}^{m-1}\sum_{s=j+1}^{m-1}h_s=(m-2)H,
$$
therefore we obtain the lower bound
$$
D_m\ge 2H^2\left( \sum_{k=1}^{2m}l_k+(m-\frac32)H \right).
$$

Thus, setting $\mathcal{L}=\sum_{k=1}^{2m}l_k$, we have the two inequalities%
$$
D_{m} \ge 2H^{2}\left(\mathcal{L}+\left( m-\frac32\right) H\right) , \quad
N_{m} \le c\left( m\right) H^{2}\left(\mathcal{L}+ mH\right) ,
$$
from which we deduce%
$$
\frac{N_m}{D_m}\le \frac12 c(m)\frac{\mathcal{L}+ mH}{\mathcal{L}+\left( m-\frac32\right) H}.
$$
For $m\ge 2$ even, it is easy to verify that $m(m-1)\le (m+2)(m-3/2)$, whence%
$$
\frac{\mathcal{L}+ mH}{\mathcal{L}+\left( m-\frac32\right) H}\le \frac{m+2}{m-1}.
$$
For $m\ge 3$ odd, one can verify that $m(m-1)\le(m+1)(m-3/2)$, whence
$$
\frac{\mathcal{L}+ mH}{\mathcal{L}+\left( m-\frac32\right) H}\le\frac{m+1}{m-1}.
$$
Finally, as $\Vert G_p^*\Vert$ is bounded above by
\begin{equation*}
\left\vert a_{m}^{\ast }\left( -p\right) \right\vert +\left\vert a_{m}^{\ast
}\left( 0\right) \right\vert +\left\vert a_{m}^{\ast }\left( p\right)
\right\vert =1+2\left( \left\vert a_{m}^{\ast }\left( -p\right)
\right\vert +\left\vert a_{m}^{\ast }\left( p\right) \right\vert \right)\\
 \leq 1+\frac{2}{m-1}\frac{N_m}{D_m},
\end{equation*}%
(the same upper bound is valid for $\left\vert
a_{i}^{\ast }\left( -p\right) \right\vert +\left\vert a_{i}^{\ast }\left(
0\right) \right\vert +\left\vert a_{i}^{\ast }\left( p\right) \right\vert $
, for all $i\in \mathbb{Z}$), we obtain respectively
$$
\Vert G_p^*\Vert\le 1+\frac{(m+2)c(m)}{\left( m-1\right) ^{2}}\;\;{\rm for}\; m\; {\rm even}
$$
$$
\Vert G_p^*\Vert\le 1+\frac{(m+1)c(m)}{\left( m-1\right) ^{2}}\;\;{\rm for}\; m\; {\rm odd}
$$
As $c(m)=\frac14m^2$ for $m$ even and $c(m)=\frac14(m^2-1)$ for $m$ odd, we obtain
$$
\Vert G_p^*\Vert\le 1+\frac14 C(m)
$$
with $C(m)=\displaystyle\frac{m^2(m+2)}{(m-1)^2}$ for $m$ even and
$C(m)= \displaystyle\frac{(m+1)^2}{m-1}$ for $m$ odd, which proves the theorem.
\quad $\square$


\section{The family $G_{p}^*$  of  integral quasi-interpolants is near-best}

We follow the notations and techniques used in section~5 for discrete QIs.
As the linear system satisfied by the coefficients of $\lambda_{i}\left(
f\right) $%
\begin{equation*}
\sum_{s=-p}^{p}\mu_{i+s}^{\left( r\right) }a_{i}\left( s\right)
=\theta_{i}^{\left( r\right) },\ 0\leq r\leq2,
\end{equation*}
is of maximal rank, it can be written as follows%
\begin{align*}
a_{i}\left( -p\right) +a_{i}\left( 0\right) +a_{i}\left( p\right) &
=1-\sum_{r\in K}a_{i}\left( r\right) , \\
\theta_{i-p}a_{i}\left( -p\right) +\theta_{i}a_{i}\left( 0\right)
+\theta_{i+p}a_{i}\left( p\right) & =\theta_{i}-\sum_{r\in K}\theta
_{i+r}a_{i}\left( r\right) , \\
\mu_{i-p}^{\left( 2\right) }a_{i}\left( -p\right) +\mu_{i}^{\left( 2\right)
}a_{i}\left( 0\right) +\mu_{i+p}^{\left( 2\right) }a_{i}\left( p\right) &
=\theta_{i}^{\left( 2\right) }-\sum_{r\in K}^{2}\mu _{i+r}^{\left( 2\right)
}a_{i}\left( r\right) ,
\end{align*}

\noindent and its general solution is%
\begin{align*}
a_{i}\left( -p\right) & =a_{i}^{\ast}\left( -p\right) -\sum_{r\in
K_{1}}\alpha_{i,r}a_{i}\left( r\right) +\sum_{s\in
K_{2}}\alpha_{i,s}a_{i}\left( s\right) , \\
a_{i}\left( 0\right) & =a_{i}^{\ast}\left( 0\right) -\sum_{r\in
K_{1}}\beta_{i,r}a_{i}\left( r\right) -\sum_{s\in
K_{2}}\beta_{i,s}a_{i}\left( s\right) , \\
a_{i}\left( p\right) & =a_{i}^{\ast}\left( p\right) +\sum_{r\in
K_{1}}\gamma_{i,r}a_{i}\left( r\right) -\sum_{s\in
K_{2}}\gamma_{i,s}a_{i}\left( s\right) ,
\end{align*}
where $a_{i}^{\ast}\left( -p\right) $, $a_{i}^{\ast}\left( 0\right) $, and $%
a_{i}^{\ast}\left( p\right) $ are the coefficients of $G_{p}^{\ast}$ studied
in section~9 above, and the various coefficients are quotiens of
determinants. Specifically,%
\begin{align*}
\alpha_{r} & =W\left( r,0,p\right) /W,\ \beta_{r}=W\left( -p,r,p\right) /W,\
\gamma_{r}=W\left( -p,r,0\right) /W, \\
\alpha_{s} & =W\left( 0,s,p\right) /W,\ \beta_{s}=W\left( -p,s,p\right) /W,\
\gamma_{s}=W\left( -p,0,s\right) /W,
\end{align*}
where $W\left( k,l,m\right) $, $k<l<m$, is the determinant with columns $%
\left( 1,\theta_{i+k},\mu_{i+k}^{\left( 2\right) }\right) ^{T}$, $\left(
1,\theta_{i+l},\mu_{i+l}^{\left( 2\right) }\right) ^{T}$, and $\left(
1,\theta_{i+m},\mu_{i+m}^{\left( 2\right) }\right) ^{T}$, and $W=W\left(
-p,0,p\right) $. As in section~5, we can write the minimization problem~P$_{%
\text{i}}$ described in section~8 (with $q=2$) in Watson's form, so we have
to minimize $\left\Vert a_{i}\right\Vert _{1}=\left\Vert a_{i}^{\ast}-A_{i}%
\widetilde{a}_{i}\right\Vert _{1}$, where we have used the same notations as
in section~5.

\noindent According to theorem 3, we must find a vector
$v^{\ast }\in \mathbb{R}^{2p+1}$ satisfying%
\begin{equation*}
\left\Vert v^{\ast }\right\Vert _{\infty }\leq 1,\quad A_{i}^{T}v^{\ast
}=0,\quad v^{\ast }\left( r\right) ={\rm sgn}\left( a_{i}^{\ast }\left(
r\right) \right) \text{ for }r=-p,0,p\text{.}
\end{equation*}%
Let us choose $v^{\ast }\left( -p\right) =-1$, $v^{\ast }\left( 0\right) =1$%
, $v^{\ast }\left( p\right) =-1$, and%
\begin{equation*}
v^{\ast }\left( j\right) =\left\{
\begin{array}{ll}
-\alpha _{i,r}+\beta _{i,r}+\gamma _{i,r}, & \text{if }j\in K_{1}, \\
\alpha _{i,r}+\beta _{i,r}-\gamma _{i,r}, & \text{if }j\in K_{2}.%
\end{array}%
\right.
\end{equation*}%
Equations $A_{i}^{T}v^{\ast }=0$ are satisfied. Moreover, as proved in
theorem~6, $\mu _{i}^{\left( 2\right) }-\theta _{i}^{\left( 2\right) }>0$
holds, and the explicit expressions for $a_{i}^{\ast }\left( -p\right) $, $%
a_{i}^{\ast }\left( 0\right) $, and $a_{i}^{\ast }\left( p\right) $, similar
to those obtained for $i=m$ in the proof of the theorem~8, imply that 
$v^{\ast}\left( r\right) ={\rm sgn}\left( a_{i}^{\ast}\left( r\right)
\right) $. 
It only remains to prove that, for $\left( r,s\right) \in
K_{1}\times K_{2}$%
\begin{equation*}
\left\vert v^{\ast }\left( r\right) \right\vert =\left\vert -\alpha
_{i,r}+\beta _{i,r}+\gamma _{i,r}\right\vert \leq 1,\quad \left\vert v^{\ast
}\left( s\right) \right\vert =\left\vert \alpha _{i,s}+\beta _{i,s}-\gamma
_{i,s}\right\vert \leq 1.
\end{equation*}%
As $\beta _{i,r}=1-\alpha _{i,r}+\gamma _{i,r}$ for $r\in K_{1}$ and $\beta
_{i,s}=1+\alpha _{i,s}-\gamma _{i,s}$ for $s\in K_{2}$, it is equivalent to
prove%
\begin{equation*}
0\leq \alpha _{i,r}-\gamma _{i,r}\leq 1,\quad 0\leq \gamma _{i,s}-\alpha
_{i,s}\leq 1,\quad \text{for }\left( r,s\right) \in K_{1}\times K_{2}\text{.}
\end{equation*}%
We only detail the proof for $i=m$ and for $r\in K_{1}$, that for $s\in K_{2}$ being quite
similar.

\noindent We prove that $\gamma _{m,r}\leq \alpha _{m,r}$, $r\in K_{1}$ by
stating that $\gamma _{m,-r}\leq \alpha _{m,-r}$, for $r\in \left\{ 1,\ldots
,p-1\right\} $, the latter being equivalent to $W\left( -p,-r,0\right) \leq
W\left( -r,0,p\right) $. By expanding the various determinants involved, we
obtain the inequality%
\begin{equation*}
\frac{\mu _{m}^{\left( 2\right) }-\mu _{m-r}^{\left( 2\right) }}{\theta
_{m}-\theta _{m-r}}\leq \frac{\mu _{m+p}^{\left( 2\right) }-\mu
_{m-p}^{\left( 2\right) }}{\theta _{m+p}-\theta _{m-p}}.
\end{equation*}%

For $1\leq i\leq m$, let $w_{i}^{\prime }=\sum_{j=i+1-r}^{i}h_{j},\;w^{\prime }=\sum_{j=1}^{m}w_{i}^{\prime},\;\tau _{i}^{\prime }=\frac{1}{m+1}\left(
\sum_{j=i-r}^{m-r}t_{j}+\sum_{j=1}^{i}t_{j}\right)$,
$w_{i}=\sum_{j=i+1-p}^{i+p}h_{j},\;\;w=\sum_{j=1}^{m}w_{i}$, and 
$\tau _{i}=\frac{1}{m+1}\left( \sum_{j=i-p}^{m-p}t_{j}+\sum_{j=1+p}^{i+p}t_{j}\right)$. \\

With these notations, we can write succesively%
\begin{align*}
\theta _{m}-\theta _{m-r}& =\frac{1}{m}\sum_{i=1}^{m}\left(
t_{i}-t_{i-r}\right) =\frac{1}{m}w^{\prime }, \\
\theta _{m+p}-\theta _{m-p}& =\frac{1}{m}\sum_{i=1}^{m}\left(
t_{i+p}-t_{i-p}\right) =\frac{1}{m}w, \\
\mu _{m}^{\left( 2\right) }-\mu _{m-r}^{\left( 2\right) }& =\frac{2}{m\left(
m+1\right) }\sum_{i=1}^{m}w_{i}^{\prime }\tau _{i}^{\prime }, \\
\mu _{m+p}^{\left( 2\right) }-\mu _{m-p}^{\left( 2\right) }& =\frac{2}{%
m\left( m+1\right) }\sum_{i=1}^{m}w_{i}\tau _{i},
\end{align*}%
and the inequality $\gamma _{m,-r}\leq \alpha _{m,-r}$ is equivalent to%
\begin{equation*}
\frac{1}{w^{\prime }}\sum_{i=1}^{m}w_{i}^{\prime }\tau _{i}^{\prime }\leq
\frac{1}{w^{\prime }}\sum_{i=1}^{m}w_{i}^{\prime }\tau _{i}^{\prime }.
\end{equation*}%
It can be interpreted as follows: the barycenter of the $m$ points $\tau
_{i}^{\prime }$ with weights $w_{i}^{\prime }$ is less than or equal to the
barycenter of the $m$ points $\tau _{i}$ with weights $w_{i}$.

\noindent The inequality $\alpha _{m,-r}-\gamma _{m,-r}\leq 1$, $r\in
\left\{ 1,\ldots ,p-1\right\} $ is equivalent to $W\left( -r,0,p\right)
-W\left( -p,-r,0\right) \leq W$ and can be written%
\begin{equation*}
\frac{\mu _{m-r}^{\left( 2\right) }-\mu _{m-p}^{\left( 2\right) }}{\theta
_{m-r}-\theta _{m-p}}\leq \frac{\mu _{m+p}^{\left( 2\right) }-\mu
_{m-r}^{\left( 2\right) }}{\theta _{m+p}-\theta _{m-r}}.
\end{equation*}%
Using the same techniques as above, we obtain successively%
\begin{align*}
\theta _{m-r}-\theta _{m-p}& =\frac{1}{m}\sum_{i=1}^{m}\left(
t_{i-r}-t_{i-p}\right) =\frac{1}{m}\overline{w}^{\prime }, \\
\theta _{m+p}-\theta _{m-r}& =\frac{1}{m}\sum_{i=1}^{m}\left(
t_{i+p}-t_{i-r}\right) =\frac{1}{m}\overline{w},
\end{align*}%
where%
\begin{equation*}
\overline{w}^{\prime }=\sum_{i=1}^{m}\overline{w}_{i}^{\prime },\quad \text{%
and}\quad \overline{w}=\sum_{i=1}^{m}\overline{w}_{i}
\end{equation*}%
with%
\begin{equation*}
\overline{w}_{i}^{\prime }=\sum_{j=i+1-p}^{i+r}h_{j},\quad \text{and}\quad
\overline{w}_{i}=\sum_{j=i+1-r}^{i+p}h_{j},
\end{equation*}%
for $1\leq i\leq m$. Defining $\overline{\tau }_{i}^{\prime }=\frac{1}{m+1}%
\left( \sum_{j=i-p}^{m-p}t_{j}+\sum_{j=1-r}^{i-r}t_{j}\right) $, and $%
\overline{\tau }_{i}^{\prime }=\frac{1}{m+1}\left(
\sum_{j=i-r}^{m-r}t_{j}+\sum_{j=1+p}^{i+p}t_{j}\right) $ for $1\leq i\leq m$%
, we get
\begin{align*}
\mu _{m-r}^{\left( 2\right) }-\mu _{m-p}^{\left( 2\right) }& =\frac{2}{%
m\left( m+1\right) }\sum_{i=1}^{m}\overline{w}_{i}^{\prime }\overline{\tau }%
_{i}^{\prime }, \\
\mu _{m+p}^{\left( 2\right) }-\mu _{m-r}^{\left( 2\right) }& =\frac{2}{%
m\left( m+1\right) }\sum_{i=1}^{m}\overline{w}_{i}\overline{\tau }_{i}.
\end{align*}%
Then the inequality $\alpha _{m,-r}-\gamma _{m,-r}\leq 1$ can be interpreted
in a geometric form as follows: the barycenter of the $m$ points $\overline{%
\tau }_{i}^{\prime }$ with weights $\overline{w}_{i}^{\prime }$ is less than
or equal to the barycenter of the $m$ points $\overline{\tau }_{i}$ with
weights $\overline{w}_{i}$.

\noindent In a similar way, one can prove that inequalities $0\leq \gamma
_{m,s}-\alpha _{m,s}\leq 1$, for all $1\leq s\leq p-1$ are equivalent to the
following inequalities:%
\begin{equation*}
\frac{\mu _{m+p}^{\left( 2\right) }-\mu _{m-p}^{\left( 2\right) }}{\theta
_{m+p}-\theta _{m-p}}\leq \frac{\mu _{m+s}^{\left( 2\right) }-\mu
_{m}^{\left( 2\right) }}{\theta _{m+s}-\theta _{m}},\quad \frac{\mu
_{m+s}^{\left( 2\right) }-\mu _{m-p}^{\left( 2\right) }}{\theta
_{m+s}-\theta _{m-p}}\leq \frac{\mu _{m+p}^{\left( 2\right) }-\mu
_{m+s}^{\left( 2\right) }}{\theta _{m+p}-\theta _{m+s}},
\end{equation*}%
and can also be interpreted in terms of barycenters of knots. Here we need
the weighted points $\left( \tau _{i}^{\prime \prime },w_{i}^{\prime \prime
}\right) $ and $\left( \overline{\tau }_{i}^{\prime \prime },\overline{w}%
_{i}^{\prime \prime }\right) $, where $\tau _{i}^{\prime \prime }=\frac{1}{%
m+1}\left( \sum_{j=i}^{m}t_{j}+\sum_{j=1+s}^{i+s}t_{j}\right) $, $\overline{%
\tau }_{i}^{\prime \prime }=\frac{1}{m+1}\left(
\sum_{j=i+s}^{m+s}t_{j}+\sum_{j=1+p}^{i+p}t_{j}\right) $, $w_{i}^{\prime
\prime }=\sum_{j=i+1+s}^{i+p}h_{j}$, and $\overline{w}_{i}^{\prime \prime
}=\sum_{j=i+1}^{i+s}h_{j}$.

\begin{theorem}
Assume that the sequence of knots $T$ satisfies, for all $i\in\mathbb{Z}$,
the following properties:

\begin{enumerate}
\item for all $1\leq r\leq p-1$, the barycenter of the $m$ points $\left(
\tau _{i}^{\prime },w_{i}^{\prime }\right) $ (resp. $\left( \overline{\tau }%
_{i}^{\prime },\overline{w}_{i}^{\prime }\right) $) is less than or equal to
the barycenter of the $m$ points $\left( \tau _{i},w_{i}\right) $ (resp. $%
\left( \overline{\tau }_{i},\overline{w}_{i}\right) $).

\item for all $1\leq s\leq p-1$, the barycenter of the $m$ points $\left(
\tau _{i},w_{i}\right) $ (resp. $\left( \overline{\tau }_{i},\overline{w}%
_{i}\right) $) is less than or equal to the barycenter of the $m$ points $%
\left( \tau _{i}^{\prime \prime },w_{i}^{\prime \prime }\right) $ (resp. $%
\left( \overline{\tau }_{i}^{\prime \prime },\overline{w}_{i}^{\prime \prime
}\right) $).
\end{enumerate}

Then, for all $i\in \mathbb{Z}$, $a_{i}^{\ast }$ is an optimal solution of
the local minimization problem~$(\bar M_i)$. Thus, for all $p\geq m$, the
spline iQIs $G_{p}^{\ast }$ of theorem~9 are near-best.
\end{theorem}

\begin{remark}
Even if the partition $T$ does not satisfy the hypothesis of theorem~9, the
operator $G_{p}^{\ast }$ is still a good iQI because its infinity norm is
uniformly bounded (theorem~8).
\end{remark}


\section{Quasi-Interpolants exact on $\PP_n$ with $n\ge 3$}

While we succeeded above in characterizing some families of near-best QIs exact on 
$\PP_2$  whose norms are {\sl uniformly bounded independently of the partition}, it is surprisingly difficult to find QIs exact on $\PP_3$ having the same property. Actually, we did not find any example of such a QI. 

For example, let us consider  the following  cubic spline dQI (which also appears in \cite{r}), defined on $I=\RR$ endowed with a non-uniform partition, by
$$
Q_3f(x)=\sum_{i\in\ZZ} \lambda_i(f) B_i(x),
$$
where $B_i$ is the cubic B-spline with support $[t_{i-2},t_{i+2}]$ centered at ${t_i}$ and
$$
\lambda_i(f)=a_i f_{i-1}+b_if_i+c_if_{i+1},
$$
where $f_i=f(t_i)$, and the coefficients are given by
$$
a_i=-\frac13\frac{h_{i}^2}{h_{i-1}(h_{i-1}+h_i)},\;\;b_i=\frac13\frac{(h_{i-1}+h_i)^2}{h_{i-1}h_i},\;\;c_i=-\frac13\frac{h_{i-1}^2}{h_{i}(h_{i-1}+h_i)}.
$$
It is easy to verify that $Q_3$ is {\sl exact on} $\PP_3$, i.e. that $a_i,b_i,c_i$ satisfy the system of linear equations:
$$
a_i+b_i+c_i=1,\;\;\;\;
t_{i-1} a_i+t_{i} b_i+t_{i+1}c_i=\theta_i,
$$
$$
t_{i-1}^2 a_i+t_{i}^2 b_i+t_{i+1}^2 c_i=\theta_i^{(2)},\;\;\;\;
t_{i-1}^3 a_i+t_{i}^3 b_i+t_{i+1}^3 c_i=\theta_i^{(3)}.
$$
This operator can also be written in the quasi-Lagrange form
$$
Q_3f(x)=\sum_{i\in\ZZ} f_i \tilde B_i(x),
$$
where the fundamental function $\tilde B_i$,  having support $[t_{i-3},t_{i+3}]$, is defined by
$$
\tilde B_i=c_{i-1} B_{i-1}+b_i B_i+a_{i+1}B_{i+1}.
$$
The Chebyshev norm $\vert \Lambda\vert_{\infty}$ of the associated Lebesgue function $\Lambda=\sum_{i\in\ZZ}\vert \tilde B_i \vert$ satisfies:
$$
\vert \Lambda \vert_{\infty}=\Vert Q_3 \Vert_{\infty}.
$$
For $x\in I=[t_2,t_3]$, we have
$\Lambda=\sum_{i=0}^{5} \vert \tilde B_i \vert$
and we shall now construct {\sl a partition for which} $\vert \Lambda\vert_{\infty}$  {\sl is arbitrary large}.
We choose $h_3=h$ as parameter and $h_i=1$ for all $i\neq 3$. Then the fundamental functions on the interval $I$ are given by :
$$
\tilde B_0=-\frac16 B_1,\; \;\tilde B_1=\frac43 B_1-\frac16 B_2,\; \;\tilde B_2=-\frac16 B_1
+\frac43 B_2-\frac{h^2}{3(1+h)} B_3,
$$
$$
\tilde B_3=-\frac16 B_2+\frac{(1+h)^2}{3h}B_3-\frac13\frac{1}{h(1+h)}B_4,
$$
$$
 \tilde B_4=-\frac{1}{3h(1+h)} B_3+\frac13\frac{(h+1)^2}{h}B_4,\;\;
 \tilde B_5=-\frac13\frac{h^2}{(1+h)}B_4.
$$
Let us compute the value of $\Lambda$ at the midpoint $s=\frac12(t_2+t_3)$ of $I$.

Setting $\alpha_1=B_1(t_1)=B_4(t_3)$, $\alpha_2=B_2(t_1)=B_3(t_3)$, 
$\delta_2=B_2(t_3)=B_3(t_2)$ and using the algebraic properties of B-splines, we get
$$
\alpha_1=\frac{h}{3(1+h)},\;\; \alpha_2=\frac{2h^3+h^2+9}{3(h+1)(h^2-h+3)},\;\;
\delta_2=\frac{h}{(h+1)(h^2-h+3)}.
$$
(we have $\alpha_1+\alpha_2+\delta_2=1$ because the B-splines sum to one).
Now, using the values of  the central BB-coefficients $\beta_2$ and $\gamma_2$ of $B_2$ on the interval $I$:
$$
\gamma_2=\frac{1}{h^2-h+3},\;\; \beta_2=1-\gamma_2,
$$
we can compute the values of  B-splines at this point
$$
B_1(s)=\frac18\alpha_1=B_4(s),\;\;B_2(s)=B_3(s)=\frac18(\alpha_2+3\beta_2+3\gamma_2+\delta_2)=\frac18(4-\alpha_1).
$$
After some algebraic calculations, we obtain  the asymptotic behaviour of the Lebesgue function at the midpoint of $I$ :  
$$
\Lambda(s)=O(h),\;\; h\to+\infty,
$$
therefore the norm of the associated QI is unbounded.\\

However, if we now assume that there exists $r>0$ such that the partition satisfies
$$
\frac1r\le \frac{h_{i+1}}{h_i}\le r,\;\; i\in\ZZ,
$$
then  we obtain the following upper bounds:
$$
\vert a_i \vert,  \vert c_i\vert \le\frac13\frac{r^2}{(1+r)},\,\;\; \vert b_i\vert \le\frac13(1+r)^2 ,
$$
from which we deduce
$$
\Vert Q_3 \Vert_{\infty}\le N(r):=\frac13 \left( (1+r)^2+\frac{2r^2}{(1+r)}\right).
$$
For example, for $r=1,2,3,4,5$, we get the following values of the upper bound of the norm :
$$
N(1)\approx 1.66,\;\;N(2)\approx 3.89,\;\;N(3)\approx 6.83,\;\;N(4)\approx 10.47,\;\;\;N(5)\approx 14.78,
$$
which are still of reasonable size. Therefore it seems  that such QIs are interesting in practice though they are not uniformly bounded with respect to $all$ partitions.



{\bf Addresses }\\

D. Barrera, M.J. Ib\'a\~nez,
Departamento de Matem\'atica Aplicada,
Facultad de Ciencias, Universidad de Granada,
Campus de Fuentenueva,
18071 GRANADA, Spain.
{\tt dbarrera@ugr.es, mibanez@ugr.es}\\

 P. Sablonni\`ere, INSA de Rennes,
20 Avenue des Buttes de Co\"esmes,
 CS 14315, 35043 RENNES Cedex, France.
{\tt psablonn@insa-rennes.fr}\\

D. Sbibih,
D\'epartement de Math\'ematiques et Informatique,
Facult\'e des Sciences, Universit\'e Mohammed 1er,
60000 OUJDA, Marocco.\\
{\tt sbibih@sciences.univ-oujda.ac.ma}

\end{document}